\documentclass[preprint,11pt]{elsarticle}  

\usepackage{amsmath,amssymb}
\usepackage{graphicx}

\journal{arXiv}

\numberwithin{equation}{subsection}

\begin{document}

\begin{frontmatter}

\title{Quadrature rules for $C^{0}$, $C^{1}$ splines,\\the real line, and the five (5) families}
\author{Helmut Ruhland}
\ead{Helmut.Ruhland50@web.de}

\begin{abstract}

The five (5) families of quadrature rules with periods of one or two intervals for the real line and spline classes $C^{0}$, $C^{1}$ are presented. The formulae allow one to calculate the points and weights of these quadrature rules in a very simple manner as for the classical Gauss-Legendre rules.

\end{abstract}


\end{frontmatter}


\section{Introduction}

In \cite{BC2016}, \cite{BC2015} and \cite{HRS2010} some quadrature rules for the real line and low degrees are presented. These rules have a period of one or two intervals (the so called ``1/2 rules'' with two different numbers of points in the two intervals).\\
In this paper, the quadrature rules for the real line and all degrees D, are derived from the explicit formulae given in \cite{RU}. First, the fixed points of the so called recursion maps are determined. Second, these fixed points are used to derive the polynomials (as sum of Gegenbauer polynomials), whose roots define the points, and to calculate the assigned weights.\\
As expected for each of the cases, continuity class $C^{0}$/$C^{1}$ of even/odd degree exists with 1 quadrature rule (altogether 4 families). The surprising exception is class $C^{1}$, of odd degree, with an additional second quadrature rule. This additional quadrature rule has slightly lower coefficients in the error term than the family with a point at the end of the interval. As the degree $D \rightarrow{} \infty{}$, these two quadrature rules converge to a common rule.

\section{$C^{0}$ quadrature rules for the real line}

Let $C_{n} (x)$ be the following Gegenbauer orthogonal polynomials for the weight function
$(1 - x^{2})$.
\begin{equation}
C_{{n}} \left( x \right) = C_{{n}}^{(3/2)} \left( x \right)
\end{equation}

\subsection{Odd degree $D = 2n - 1$, periodic with two intervals, the ``1/2 rules''}

The rule has a period of 2 intervals, define the first interval by:
\begin{align}
R_{{n}} \left( x \right) & = {n}^{2}C_{{n}} \left( x \right) - \left( n+1
 \right) ^{2}C_{{n-2}} \left( x \right) \\
S_{{n-1}} \left( x \right) & = n\,C_{{n-1}} \left( x \right) - \left( n+1 \right) x\,C_{{n-2}} \left( x \right) \\
A & = 2\, \left( n+1 \right) \left( 2\,n+1 \right) {n}^{2}
\end{align}
Let the $n$ points of this interval be:
\begin{equation}
x_{1}, \dots, x_{n} = \mbox{the roots of the polynomial $R_{n} (x)$}
\end{equation}
Let the $n$ weights be:
\begin{equation}
w_{{i}}={\frac {A} {{R_{{n}}}' \left( x_{{i}} \right) S_{{n-1}} \left( x_{{i}} \right) }}
\end{equation}

Define the second interval by:
\begin{align}
R_{{n-1}} \left( x \right) & = C_{{n-1}} \left( x \right) \\
S_{{n-2}} \left( x \right) & = \left( 2\,n-1 \right) C_{{n-2}} \left( x \right) -n\,x\,C_{{n-3}} \left( x \right) \\
A & = 2\, n \left( 2\,n-1 \right)
\end{align}
Let the $n - 1$ points of this interval be:
\begin{equation}
x_{1}, \dots, x_{n-1} = \mbox{the roots of the polynomial $R_{n-1} (x)$}
\end{equation}
Let the $n - 1$ weights be:
\begin{equation}
w_{{i}}={\frac {A}{ {R_{{n-1}}}'  \left( x_{{i}} \right) S_{{n-2}} \left( x_{{i}} \right) }}
\end{equation}

For low $n$, these are the well-known $S_{2n-1,0}$ quadrature rules for the real line 
(see \cite[pp. 18-21]{BC2015} for the cases $n$ = 3, 4, 5).

\subsection{Even degree $D = 2n$, periodic with one  interval}

The rule has a period of 1 interval, define this interval by:
\begin{align}
\delta & = \sqrt {{\frac {n+2}{n}}} \label{deltabothsigns} \\
R_{{n}} \left( x \right) & = C_{{n}} \left( x \right) +\delta\,C_{{n-1}}
 \left( x \right) \\
S_{{n}} \left( x \right) & = \left( 2\,n+1+\delta\,x\,n \right) C_{{n-1}}
 \left( x \right) - \left( n+1 \right) x\,C_{{n-2}} \left( x \right) \\
A & = 2\, \left( n+1 \right) \left( 2\,n+1 \right)
\end{align}
Let the $n$ points of this interval be:
\begin{equation}
x_{1}, \dots, x_{n} = \mbox{the roots of the polynomial $R_{n} (x)$}
\end{equation}
Let the $n$ weights be:
\begin{equation}
w_{{i}}={\frac {A}{ {R_{{n}}}' \left( x_{{i}} \right) S_{{n}} \left( x_{{i}} \right) }}
\end{equation}

Note that in the formula for $\delta$ in (\ref{deltabothsigns}) both signs of the square root can be choosen, this rule has  no reflection symmetry and changing the sign is equivalent to a reflection.

\section{$C^{1}$ quadrature rules for the real line}

Let $C{}_{n} (x)$ be the following Gegenbauer orthogonal polynomials for the weight function
$(1 - x^{2})^{2}$
\begin{equation}
C_{{n}} \left( x \right) = C_{{n}}^{(5/2)} \left( x \right)
\end{equation}

\subsection{Odd degree $D = 2n + 1$, periodic with one interval}

For this case, there exist two quadrature rules.

\subsubsection{The first quadrature rule with a point at the end of an
interval}

The rule has a period of 1 interval, define this interval by:
\begin{align}
R_{{n-1}} \left( x \right) & = C_{{n-1}} \left( x \right) \\
S_{{n-2}} \left( x \right) & = C_{{n-2}} \left( x \right) \\
A & = 2\, n \left( n+1 \right)  \left( n+2 \right) / \, 9
\end{align}
Let the $n$ points of this interval be:
\begin{align} \begin{split}
x_{1} & = -1 \\
x_{2}, \dots, x_{n} & = \text{the roots of the polynomial $R_{n-1} (x)$}
\end{split}
\label{easypoints}
\end{align}
Let the $n$ weights be:
\begin{align}
w_{{1}} & = {\frac {16\, \left( 2\,{n}^{2}+6\,n+1 \right)}{3\,n \left( n+1 \right)  \left( 
n+2 \right)  \left( n+3 \right) }}
\label{noteasyweight} \\
w_{{i}} & = {\frac {A}{ {R_{{n-1}}}'  \left( x_{{i}} \right) S_{{n-2}} \left( x_{{i}} \right) \left( 1-{x_{{i}}}^{2} \right) ^{2}}}  \qquad i = 2, \dots, n \\
        & = {\frac {2\,n \left( n+1 \right)  \left( n+2 \right)}{9\, {C_{{n-1}}}'  \left( x_{{i}} \right) C_{{n-2}} \left( x_{{i}} \right) \left( 1-{x_{{i}}}^{2} \right) ^{2}}}
\label{easyweights}
\end{align}

For low $n$, these are the well-known S$_{2n+1,1}$ quadrature rules for the real line 
(see \cite[pp. 13-17]{BC2015} for the cases $n$ = 2, 3, 4). \bigskip

\noindent
Remark to a proof of this rule:\\
Using the Fundamental Theorem (FT) for Gaussian quadrature with the weight function
$\omega = (1 - x^{2})^{2}$ and $n - 1$ points we get (\ref{easypoints}) and
(\ref{easyweights}), when we insert the additonal factor $(1 - x_i^{2})^{2}$ in the
denominator of the weight. This additional factor is is necesary because the FT states
for the weighted integral $\int \omega (x) h (x) dx = \sum w_i h (x_i)$. But we want to
interpret $\omega (x) h (x)$ as spline functon with a support of 1 interval and want
to evaluate in the sum at this spline function instead of an evaluation of $h (x_i)$.\\
These $n-1$ points/weights allow the exact calculation for the splines with
a support of 1 interval up to the degree:\\
$n-1$ (the degree of $C_{n-1}$) + $n-2$ (all $C_{0}$ ... $C_{n-2}$ are orthogonal to 
$C_{n-1}$) + $4$ (the degree of the weight function) = $2n+1$ (the neccesary degree
for a Gaussian rule).\\
The additional point at -1 does not change anything for the support 1 splines, they are 0 at x = -1.
Because these $n-1$ points
are reflection symmetric the quadrature for the odd spline with support 2 is also correct. The weight $w_1$
can now be determined so that the quadrature for the even spline with support 2 is correct too.
But for me the value (\ref{noteasyweight}) of $w_1$ is far from being obvious (if you have a simple 
proof, please let me know).

\subsubsection{The second quadrature rule}

The rule has a period of 1 interval, define this interval by:
\begin{align}
\delta & = \sqrt {{\frac {3\,\left( {n}^{2}+3\,n-1 \right)}{n \left( n+3 \right) }}}
\label{deltaplussign1} \\
\begin{split}
R_{{n}} \left( x \right)  & = \left( n-1 \right) \left( 2\,{n}^{2}+2\,n-3 \right) C_{{n}} \left( x \right) \\
 & \quad - \left( n+3 \right) \left( 2\,{n}^{2}+6\,n+7-2\, \left( 2\,n+3 \right) \delta \right) C_{{n-2}} \left( x \right)
\end{split} \\
\begin{split}
S_{{n+1}} \left( x \right) & = n \left( 6\,{n}^{2}+6\,n-3+2\, \left( 2\,n+1 \right) \delta \right) \left( 1+{x}^{2} \right) C_{{n-1}} \left( x \right) \\
 & \quad -4 \, \left( 2\,n+1 \right) \left( 2\,{n}^{2}+6\,n+1 \right) x \, C_{{n-2}} \left( x \right) \\
 & \quad + \left( n+2 \right)  \left( 2\,{n}^{2}+6\,n+1 \right)  \left( 1+{x}^{2} \right) C_{{n-3}} \left( x \right)
\end{split} \\
\begin{split}
A & = 2 \, \left( n-1 \right) \left( n+1 \right)  \left( n+2 \right) \left( 2\,n+1 \right) \left( 2\,n+3 \right) \\
 & \quad \cdot  \left( 2\,{n}^{2}+2\,n-3 \right)  \left( 2\,{n}^{2}+6\,n+1 \right) / \, 9
\end{split}
\end{align}

Let the $n$ points of this interval be:
\begin{equation}
x_{1}, \dots, x_{n} = \mbox{the roots of the polynomial $R_{n} (x)$}
\end{equation}
Let the $n$ weights be:
\begin{eqnarray}
w_{{i}} & = & {\frac {A}{ {R_{{n}}}'  \left( x_{{i}} \right) S_{{n+1}} \left( x_{{i}} \right) }}
\end{eqnarray}

Note that in the formulae for $\delta$ in (\ref{deltaplussign1}) and the following (\ref{deltaplussign2}) only the positive square root is to be choosen, because using the negative square root results in some roots of $R_{n} (x)$ falling outside the interval [-1, +1] and so results in no quadrature rule.
 
\subsection{Even degree $D = 2n$, periodic with two intervals - the ``1/2 rules''}

The rule has a period of 2 intervals, define the first interval by:
\begin{align}
\delta & =  \sqrt {3\,n \left( n+2 \right)  \left( {n}^{2}+2\,n-2  \right) }
\label{deltaplussign2} \\
\begin{split}
R_{{n-1}} \left( x \right) & = \left( n-1 \right) \left( 2\,{n}^{2}+2\,n-3 \right) C_{{n-1}} \left( x \right) \\
 & \quad + \left( 2\,\delta+3-n-6\,{n}^{2}-2\,{n}^{3} \right) C_{{n-2}} \left( x \right)
\end{split} \\
\begin{split}
S_{{n-2}} \left( x \right) & = \left( 3\, \left( n+2 \right)  \left( 2\,
{n}^{2}-1 \right) -2\,\delta \right) C_{{n-2}} \left( x \right) \\
 & \quad +  \left( n+2 \right)  \left( 2\,{n}^{2}+2\,n-3 \right) C_{{n-3}}
 \left( x \right)
\end{split} \\
A & =  2 \, \left( n-1 \right) n \left( n+1 \right) \left( n+2 \right) \left( 2\,n+1 \right) \left( 2\,{n}^{2}+2\, n-3 \right) ^{2} / \, 9
\end{align}
Let the $n$ points of this first interval be:
\begin{align} \begin{split}
x_{1} & = -1 \\
x_{2}, \dots, x_{n} & = \text{the roots of the polynomial $R_{n-1} (x)$}
\end{split} \end{align}
Let the $n$ weights be:
\begin{align}
\begin{split}
w_{{1}} & = {\frac { 8\, \left( 2\,{n}^{2}+4\,n-3 \right) \left( 2\,{n}^{4}+8\,{n}^{3}+4\,{n}^{2}-8\,n-3-\delta \right) }
{ 3\, \left( n-1 \right) n \left( n+2 \right) \left( n+3 \right)  \left( {n}^{2}+2\,n-2
 \right)  \left( n+1 \right) ^{2} }} 
\end{split} \\
w_{{i}} & = {\frac {A}{ {R_{{n-1}}}'  \left( x_{{i}} \right) S_{{n-2}} \left( x_{{i}} \right) \left( 1+x_{{i}} \right)  \left( 1-x_{{i}} \right) ^{2}}} \qquad i = 2, \dots, n 
\end{align}

Define the second interval by:\medskip

This interval has $n - 1$ points, which are the reflection at 0, i.e $x_{i} \rightarrow{} - x_{i}$ of the points $x_{2}, \dots, x_{n}$ of the first interval.\\
$( \, \dots \, )_{1.}$ in the following formulae on the rhs means: take in the bracket the points/weights of the first interval\\

Let the $n - 1$ points of this interval be:
\begin{equation}
x_{1}, \dots, x_{n-1} = ( \, - \, x_{2}, \dots, - \, x_{n} \, )_{1.}
 \text{\quad 1. interval, reflected at 0}
\end{equation}

Let the $n - 1$ weights of this interval be:
\begin{equation}
w_{1}, \dots, w_{n-1} = ( \, w_{2}, \dots, w_{n} \, )_{1.}
\end{equation}

For low $n$, these are the well-known S$_{2n,1}$ quadrature rules for the real line 
(see \cite[p. 308]{HRS2010} for the ``two-third'' quartic case $n$ = 2 and \cite[p. 19]{BC2016} for the ``two-and-half'' sextic case $n$ = 3).

\bigskip
\bigskip
This research did not receive any specific grant from funding agencies in the public, commercial, or not-for-profit sectors.


\bibliographystyle{elsarticle-num}

\appendix

\section{The rules as code for the Computer Algebra System MAPLE}

A quadrature rule is displayed in the form: [ [ $x_{1}$, $w_{1}$ ], [ $x_{2}$, $w_{2}$ ], ... ], a list of lists in the langauge of CAS`s.\\
If the rules can be presented by (nested) square roots, the algebraic numbers are given. Otherwise, we give 25 significant digits of the float representation.\\

For rules with a period of one interval the points/weights, calculated with the formulae above, are scaled to the interval [0, 1]. For rules with a period of two intervals the points/weights are scaled to [0, 1] and [1, 2].\\
When implementing the above formulae do not forget to scale the weights with $interval \, length \, / \, 2$. 

\subsection{The class $C^{0}$, real line rules}

\fontsize{5pt}{6pt} \ttfamily


\medskip \noindent
C0xD2 := [\\ \relax
 [1/2-1/6*3\^{}(1/2), 1] ];\\ \relax
C0xD3 := [\\ \relax
 [-1/4*2\^{}(1/2)+1/2, 2/3], [1/4*2\^{}(1/2)+1/2, 2/3],\\ \relax
 [3/2, 2/3] ];\\ \relax
C0xD4 := [\\ \relax
 [-1/10*2\^{}(1/2)-1/10*7\^{}(1/2)+1/2, -1/84*7\^{}(1/2)*2\^{}(1/2)+1/2],\\ \relax
 [-1/10*2\^{}(1/2)+1/10*7\^{}(1/2)+1/2, 1/84*7\^{}(1/2)*2\^{}(1/2)+1/2] ];\\ \relax
C0xD5 := [ \# see Barton, Calo [2], C0 Quintics, d = 5, c = 0, formulae (40) and (41)\\ \relax
 [-1/30*165\^{}(1/2)+1/2, 15/44], [1/2, 16/33], [1/30*165\^{}(1/2)+1/2, 15/44],\\ \relax
 [3/2-1/10*5\^{}(1/2), 5/12], [3/2+1/10*5\^{}(1/2), 5/12] ];\\ \relax
C0xD6 := [\\ \relax
 [.0529116719292753211479553, .2586020762640311600074680], [.3897721580810322272553957, .4016943058349034864147465],\\ \relax
 [.7806745024034483883695585, .3397036179010653535777854] ];\\ \relax
C0xD7 := [ \# see Barton, Calo [2], C0 Septics, d = 7, c = 0, formula (42), above and below\\ \relax
 [-1/56*(378+14*393\^{}(1/2)) \^{}(1/2)+1/2, 11/40-163/47160*393\^{}(1/2)],\\ \relax
 [-1/56*(378-14*393\^{}(1/2)) \^{}(1/2)+1/2, 11/40+163/47160*393\^{}(1/2)],\\ \relax
 [1/56*(378-14*393\^{}(1/2)) \^{}(1/2)+1/2, 11/40+163/47160*393\^{}(1/2)],\\ \relax
 [1/56*(378+14*393\^{}(1/2)) \^{}(1/2)+1/2, 11/40-163/47160*393\^{}(1/2)],\\ \relax
 [3/2-1/14*21\^{}(1/2), 49/180], [3/2, 16/45], [3/2+1/14*21\^{}(1/2), 49/180]];\\ \relax
C0xD8 := [\\ \relax
 [.0338755371265535325717510, .1662405692717599599932669], [.2616435694284255998558377, .2876981623645600825790484],\\ \relax
 [.5749604515321156967840158, .3183454907229987216856932], [.8573549149369964932109196, .2277157776406812357419916] ];\\ \relax
C0xD9 := [ \# see Barton, Calo [2], C0 Nonics, d = 9, c = 0, table 6, \#el. 3 and 4\\ \relax
 [-1/210*(6615+210*231\^{}(1/2))\^{}(1/2)+1/2, 28/145-23/6380*231\^{}(1/2)],\\ \relax
 [-1/210*(6615-210*231\^{}(1/2))\^{}(1/2)+1/2, 28/145+23/6380*231\^{}(1/2)],\\ \relax
 [1/2, 128/435], \\ \relax
 [1/210*(6615-210*231\^{}(1/2))\^{}(1/2)+1/2, 28/145+23/6380*231\^{}(1/2)],\\ \relax
 [1/210*(6615+210*231\^{}(1/2))\^{}(1/2)+1/2, 28/145-23/6380*231\^{}(1/2)],\\ \relax
 [3/2-1/42*(147+42*7\^{}(1/2))\^{}(1/2), 7/30-1/60*7\^{}(1/2)],\\ \relax
 [3/2-1/42*(147-42*7\^{}(1/2))\^{}(1/2), 7/30+1/60*7\^{}(1/2)],\\ \relax
 [3/2+1/42*(147-42*7\^{}(1/2))\^{}(1/2), 7/30+1/60*7\^{}(1/2)],\\ \relax
 [3/2+1/42*(147+42*7\^{}(1/2))\^{}(1/2), 7/30-1/60*7\^{}(1/2)] ];\\ \relax
C0xD10 := [\\ \relax
 [.0235295413811732733578363, .1157259571662055027409935], [.1864032815050772945242715, .2117953691406563246913401],\\ \relax
 [.4299526855497696596273958, .2642347893930234698014511], [.6911353204896879120519700, .2461539393851279315081472],\\ \relax
 [.9000664536606729494127386, .1620899449149867712580682] ];\\ \relax
C0xD11 := [\\ \relax
 [.0202715343459102791524548, .0990651582961091991162324], [.1612217078414028759175214, .1850629079215814540702768],\\ \relax
 [.3777774454048602455845435, .2396814575918331563373002], [.6222225545951397544154565, .2396814575918331563373002],\\ \relax
 [.8387782921585971240824786, .1850629079215814540702768], [.9797284656540897208475452, .0990651582961091991162324],\\ \relax
 [1.084888051860716535063984, .1384130236807829740053502], [1.265575603264642893098114, .2158726906049313117089355],\\ \relax
 [1.500000000000000000000000, .2438095238095238095238095], [1.734424396735357106901886, .2158726906049313117089355],\\ \relax
 [1.915111948139283464936016, .1384130236807829740053502] ];\\ \relax
C0xD12 := [\\ \relax
 [.0172892396789921364872979, .0851485943344666872126749], [.1390618167952894575338889, .1610397260440442134913879],\\ \relax
 [.3300658706736885789701525, .2146144843199893151992085], [.5540853769912574828028730, .2258896080582861900223107],\\ \relax
 [.7668399370334643249823590, .1924535819133622825479249], [.9261884038167115125268984, .1208540053298513115264931] ];\\ \relax
C0xD13 := [\\ \relax
 [.0151859177411498769190688, .0744264776909320504064757], [.1223855762277330681268381, .1426255359977720749903832],\\ \relax
 [.2933787936045333156804095, .1944155187788283421356086], [.5000000000000000000000000, .2127792207792207792207792],\\ \relax
 [.7066212063954666843195905, .1944155187788283421356086], [.8776144237722669318731619, .1426255359977720749903832],\\ \relax
 [.9848140822588501230809312, .0744264776909320504064757], \\ \relax
 [1.064129925745196692331277, .1053521135717530196914960], [1.204149909283428848927745, .1705613462417521823821203],\\ \relax
 [1.395350391048760565615671, .2062293973293519407835265], [1.604649608951239434384329, .2062293973293519407835265],\\ \relax
 [1.795850090716571151072255, .1705613462417521823821203], [1.935870074254803307668723, .1053521135717530196914960] ];\\ \relax
C0xD14 := [\\ \relax
 [.0132382049859402066250271, .0652545015152810717772603], [.1075275184163180497279689, .1260332526421032697696779],\\ \relax
 [.2599422611502342475958597, .1749288046650861619891666], [.4485070043501314855627717, .1973499794883650908899665],\\ \relax
 [.6445836536136272739800257, .1897525091922582802460465], [.8183336274094846945391861, .1532689843692582894782297],\\ \relax
 [.9432925989678049829189991, .0934119681276478358496524] ]; \\ \relax
C0xD15 := [\\ \relax
 [.0118021490241755555547310, .0579536270034359693949455], [.0959578630301719500366127, .1129430902900632556391664],\\ \relax
 [.2335191673685803759775308, .1590935465188152541807088], [.4069904537399158205303826, .1838986250765744096740681],\\ \relax
 [.5930095462600841794696174, .1838986250765744096740681], [.7664808326314196240224692, .1590935465188152541807088],\\ \relax
 [.9040421369698280499633872, .1129430902900632556391664], [.9881978509758244444452690, .0579536270034359693949455],\\ \relax
 [1.050121002294269921343827, .0827476807804027625231698], [1.161406860244631123277057, .1372693562500808676403528],\\ \relax
 [1.318441268086910920644624, .1732142554865231725575658], [1.500000000000000000000000, .1857596371882086167800454],\\ \relax
 [1.681558731913089079355376, .1732142554865231725575658], [1.838593139755368876722943, .1372693562500808676403528],\\ \relax
 [1.949878997705730078656173, .0827476807804027625231698] ];\\ \relax

\normalsize \normalfont

\subsection{The class $C^{1}$, real line rules}

\fontsize{5pt}{6pt} \ttfamily

\medskip \noindent
C1xD3 := [\\ \relax
 [0, 1] ];\\ \relax
C1xD3x2 := [\\ \relax
 [1/2, 1] ];\\ \relax
C1xD4 := [ \# see Hughes, Reali, Sangalli [4], formula (29)\\ \relax
 [0, 13/20], [2/3, 27/40],\\ \relax
 [4/3, 27/40] ];\\ \relax
C1xD5 := [ \# see Barton, Calo [2], C1 Quintics, d = 5, c = 1, formula (33)\\ \relax
 [0, 7/15], [1/2, 8/15] ];\\ \relax
C1xD5x2 := [\\ \relax
 [-1/30*(225-30*30\^{}(1/2))\^{}(1/2)+1/2, 1/2], [1/2+1/30*(225-30*30\^{}(1/2))\^{}(1/2), 1/2] ];\\ \relax
C1xD6 := [ \# see Barton, Calo [1], C1 Sextics, d = 6, c = 1, formula (37) and (38)\\ \relax
 [0, 387/1040-3/1040*65\^{}(1/2)],\\ \relax
 [67/98-1/98*65\^{}(1/2)-3/98*78\^{}(1/2)+1/98*30\^{}(1/2),\\ \relax
  -673/99840*5\^{}(1/2)*6\^{}(1/2)+1693/4160+2047/299520*6\^{}(1/2)*13\^{}(1/2)+3/4160*5\^{}(1/2)*13\^{}(1/2)],\\ \relax
 [67/98-1/98*65\^{}(1/2)+3/98*78\^{}(1/2)-1/98*30\^{}(1/2),\\ \relax
   673/99840*5\^{}(1/2)*6\^{}(1/2)+1693/4160-2047/299520*6\^{}(1/2)*13\^{}(1/2)+3/4160*5\^{}(1/2)*13\^{}(1/2)],\\ \relax
\# reflection of the points (without 0) of the previous interval [0, 1] at 1\\ \relax  
 [129/98+1/98*65\^{}(1/2)-3/98*78\^{}(1/2)+1/98*30\^{}(1/2),\\ \relax
  673/99840*5\^{}(1/2)*6\^{}(1/2)+1693/4160-2047/299520*6\^{}(1/2)*13\^{}(1/2)+3/4160*5\^{}(1/2)*13\^{}(1/2)],\\ \relax
 [129/98+1/98*65\^{}(1/2)+3/98*78\^{}(1/2)-1/98*30\^{}(1/2),\\ \relax
  -673/99840*5\^{}(1/2)*6\^{}(1/2)+1693/4160+2047/299520*6\^{}(1/2)*13\^{}(1/2)+3/4160*5\^{}(1/2)*13\^{}(1/2)] ];\\ \relax
C1xD7 := [ \# see Barton, Calo [2], C1 Septics, d = 7, c = 1, formulae (35) and (37))\\ \relax
 [0, 37/135], [-1/14*7\^{}(1/2)+1/2, 49/135], [1/14*7\^{}(1/2)+1/2, 49/135] ];\\ \relax
C1xD7x2 := [\\ \relax
 [-1/14*(45-2*102\^{}(1/2)) \^{}(1/2)+1/2, 659/2310+8/3465*102\^{}(1/2)], [1/2, 496/1155-16/3465*102\^{}(1/2)],\\ \relax
 [1/2+1/14*(45-2*102\^{}(1/2))\^{}(1/2), 659/2310+8/3465*102\^{}(1/2)] ];\\ \relax
C1xD8 := [\\ \relax
 [                        0., .2193074680769989614308703], [.2530396486742266102610894, .3026115552696555127214919],\\ \relax
 [.5824768844382089931611331, .3350741306497649127026503], [.8828032942621445821952574, .2526605800420800938604227],\\ \relax
\# reflection of the points of the previous interval\\ \relax  
 [1.117196705737855417804743, .2526605800420800938604227], [1.417523115561791006838867, .3350741306497649127026503],\\ \relax
 [1.746960351325773389738911, .3026115552696555127214919] ];\\ \relax
C1xD9 := [ \# see Barton, Calo [2], C1 Nonics, d = 9, c = 1, formulae (38) and (39)\\ \relax
 [0, 19/105], [1/2-1/6*3\^{}(1/2), 9/35], [1/2, 32/105], [1/2+1/6*3\^{}(1/2), 9/35] ];\\ \relax
C1xD9x2 := [\\ \relax
 [-1/1554*(310023-32634*7\^{}(1/2)+1554*(27930-6048*7\^{}(1/2)) \^{}(1/2))\^{}(1/2)+1/2,\\ \relax
  -1/642320*(27930-6048*7\^{}(1/2))\^{}(1/2)*(249-3*7\^{}(1/2))+1/4],\\ \relax
 [-1/1554*(310023-32634*7\^{}(1/2)-1554*(27930-6048*7\^{}(1/2))\^{}(1/2))\^{}(1/2)+1/2,\\ \relax
   1/4+1/642320*(27930-6048*7\^{}(1/2))\^{}(1/2)*(249-3*7\^{}(1/2))],\\ \relax
 [1/2+1/1554*(310023-32634*7\^{}(1/2)-1554*(27930-6048*7\^{} (1/2))\^{}(1/2))\^{}(1/2),\\ \relax
  1/4+1/642320*(27930-6048*7\^{}(1/2))\^{}(1/2)*(249-3*7\^{}(1/2))],\\ \relax
 [ 1/2+1/1554*(310023-32634*7\^{}(1/2)+1554*(27930-6048*7\^{}(1/2))\^{}(1/2))\^{}(1/2),\\ \relax
  -1/642320*(27930-6048*7\^{}(1/2))\^{}(1/2)*(249-3*7\^{}(1/2))+1/4] ];\\ \relax
C1xD10 := [\\ \relax
 [                        0., .1509390484297939060474049], [.1779083893587173174142483, .2193161134540970831046591],\\ \relax
 [.4298451600684519595332682, .2732671011137680849722728], [.7001314274491368990752840, .2551024072871098616054254],\\ \relax
 [.9188068084330307451804244, .1768448539301280172939403],\\ \relax
\# reflection of the points of the previous interval\\ \relax  
 [1.081193191566969254819576, .1768448539301280172939403], [1.299868572550863100924716, .2551024072871098616054254],\\ \relax
 [1.570154839931548040466732, .2732671011137680849722728], [1.822091610641282682585752, .2193161134540970831046591] ];\\ \relax
C1xD11 := [\\ \relax
 [0, 9/70], [-1/66*(297+132*3\^{}(1/2))\^{}(1/2)+1/2, 61/280-9/560*3\^{}(1/2)],\\ \relax
 [-1/66*(297-132*3\^{}(1/2))\^{}(1/2)+1/2, 1/2*(13503687120-996173640*3\^{}(1/2))/(32020833600-4648492800*3\^{}(1/2))],\\ \relax
 [ 1/66*(297-132*3\^{}(1/2))\^{}(1/2)+1/2, 1/2*(13503687120-996173640*3\^{}(1/2))/(32020833600-4648492800*3\^{}(1/2))],\\ \relax
 [1/66*(297+132*3\^{}(1/2))\^{}(1/2)+1/2, 61/280-9/560*3\^{}(1/2)] ];\\ \relax
C1xD11x2 := [\\ \relax
 [-1/1254*(235125-3762*130\^{}(1/2)+1254*(12030-498*130\^{}(1/2))\^{}(1/2))\^{}(1/2)+1/2,\\ \relax
   3571/19404+32/121275*130\^{}(1/2)-1/63817815600*(12030-498*130\^{}(1/2))\^{}(1/2)*(30199010-96347*130\^{}(1/2))],\\ \relax
 [-1/1254*(235125-3762*130\^{}(1/2)-1254*(12030-498*130\^{}(1/2))\^{}(1/2))\^{}(1/2)+1/2,\\ \relax
   3571/19404+32/121275*130\^{}(1/2)+1/63817815600*(12030-498*130\^{}(1/2))\^{}(1/2)*(30199010-96347*130\^{}(1/2))],\\ \relax
 [1/2, 1280/4851-128/121275*130\^{}(1/2)],\\ \relax
 [1/2+1/1254*(235125-3762*130\^{}(1/2)-1254*(12030-498*130\^{}(1/2))\^{}(1/2))\^{}(1/2),\\ \relax
  3571/19404+32/121275*130\^{}(1/2)+1/63817815600*(12030-498*130\^{}(1/2)) \^{}(1/2)*(30199010-96347*130\^{}(1/2))],\\ \relax
 [1/2+1/1254*(235125-3762*130\^{}(1/2)+1254*(12030-498*130\^{}(1/2))\^{}(1/2))\^{}(1/2),\\ \relax
  3571/19404+32/121275*130\^{}(1/2)-1/63817815600*(12030-498*130\^{}(1/2))\^{}(1/2)*(30199010-96347*130\^{}(1/2))] ];\\ \relax
C1xD12 := [\\ \relax
 [                        0., .1103223678450469948505275], [.1317058754450282678835884, .1652544972251845895290208],\\ \relax
 [.3274845988167730244040532, .2198877806667871052807681], [.5570152583779801198333197, .2314870279444837967315375],\\ \relax
 [.7751968110233312377364702, .1976531719862662098657520], [.9404227556744911320335707, .1305563382547548011676578],\\ \relax
\# reflection of the points of the previous interval\\ \relax  
 [1.059577244325508867966429, .1305563382547548011676578], [1.224803188976668762263530, .1976531719862662098657520],\\ \relax
 [1.442984741622019880166680, .2314870279444837967315375], [1.672515401183226975595947, .2198877806667871052807681],\\ \relax
 [1.868294124554971732116412, .1652544972251845895290208] ];\\ \relax
C1xD13 := [\\ \relax
 [0, 109/1134], [-1/286*(7865+572*55\^{}(1/2))\^{}(1/2)+1/2, 2783/16200-803/226800*55\^{}(1/2)],\\ \relax
 [-1/286*(7865-572*55\^{}(1/2))\^{}(1/2)+1/2, 2783/16200+803/226800*55\^{}(1/2)],\\ \relax
 [1/2, 1024/4725], [1/286*(7865-572*55\^{}(1/2))\^{}(1/2)+1/2, 2783/16200+803/226800*55\^{}(1/2)],\\ \relax
 [1/286*(7865+572*55\^{}(1/2))\^{}(1/2)+1/2, 2783/16200-803/226800*55\^{}(1/2)] ];\\ \relax
C1xD13x2 := [\\ \relax
 [.0517743418529121021303818, .1137024757183161908921092], [.1968522913811618340107702, .1750639562602571514151027],\\ \relax
 [.3928206875258624010668831, .2112335680214266576927880], [.6071793124741375989331169, .2112335680214266576927880],\\ \relax
 [.8031477086188381659892298, .1750639562602571514151027], [.9482256581470878978696182, .1137024757183161908921092] ];\\ \relax
C1xD14 := [\\ \relax
 [                        0., .0841900586506389791515548], [.1013410600943745316710965, .1285896866359605905285117],\\ \relax
 [.2566656791020788952980904, .1781884212907189038041685], [.4487283186403862073586177, .2010095393585452378524983],\\ \relax
 [.6484659216760909534796831, .1933390680669159773453261], [.8256244860059208689942391, .1565039141154620584005221],\\ \relax
 [.9544202953318887983522730, .1002743412070777424931959],\\ \relax
\# reflection of the points of the previous interval\\ \relax  
 [1.045579704668111201647727, .1002743412070777424931959], [1.174375513994079131005761, .1565039141154620584005221], \\ \relax
 [1.351534078323909046520317, .1933390680669159773453261], [1.551271681359613792641382, .2010095393585452378524983], \\ \relax
 [1.743334320897921104701910, .1781884212907189038041685], [1.898658939905625468328904, .1285896866359605905285117] ]; \\ \relax
C1xD15 := [\\ \relax
 [0.                        , .0746031746031746031746031], [.0900770022682565225919989, .1147702511103050197233797],\\ \relax
 [.2296976813063206480102970, .1613996351720511671462149], [.4056612887546070310868762, .1865285264160565115431038],\\ \relax
 [.5943387112453929689131238, .1865285264160565115431038], [.7703023186936793519897030, .1613996351720511671462149],\\ \relax
 [.9099229977317434774080010, .1147702511103050197233797] ];\\ \relax
C1xD15x2 := [\\ \relax
 [.0403016742861217670180007, .0887918380491571945747060], [.1549756285834430320212064, .1401430685827084516329161],\\ \relax
 [.3150674953504518170684371, .1764672634751373029971663], [.5000000000000000000000000, .1891956597859941015904230],\\ \relax
 [.6849325046495481829315629, .1764672634751373029971663], [.8450243714165569679787936, .1401430685827084516329161],\\ \relax
 [.9596983257138782329819993, .0887918380491571945747060] ];\\ \relax

\normalsize \normalfont

\subsection{The constants $\delta{}$, A and the polynomials R, S}

\fontsize{5pt}{6pt} \ttfamily

\medskip \noindent
\#\#\#\#\#\#\#\#\#\#\#\#\#\#\#\#\#\#\#\#\#\#\#\#\#\#\#\#\#\#\#\#\#\#\#\#\#\#\#\#\#\#\#\#\\
\# two families of C0 rules, see section 2. \#\\
\#\#\#\#\#\#\#\#\#\#\#\#\#\#\#\#\#\#\#\#\#\#\#\#\#\#\#\#\#\#\#\#\#\#\#\#\#\#\#\#\#\#\#\#\\
\#\\
\# get the orthogonal Gegenbauer polynomials for the weight (1 - x \^{} 2) with\\
\# C[n] := simplify (GegenbauerC (n, 3 / 2, x));\\
\#\\
\# C0, periodic with 2 intervals, odd degree 2 * n - 1, "1/2 rules"\\
\# see section 2.1.\\
\#\\
\# 1. interval with n points\\
R[n]     = n \^{} 2 * C[n] - (n + 1) \^{} 2 * C[n - 2];\\
S[n - 1] = n * C[n - 1] - (n + 1) * x * C[n - 2];\\
A = 2 * (n + 1) * (2 * n + 1) * n \^{} 2;\\
\# 2. interval with n - 1 points\\
R[n - 1] = C[n - 1];\\
S[n - 2] = (2 * n - 1) * C[n - 2] - n * x * C[n - 3];\\
A = 2 * n * (2 * n - 1);\\
\#\\
\# C0, periodic with 1 interval, even degree 2 * n\\
\# see section 2.2.\\
\#\\
delta = sqrt ((n + 2) / n);\\
R[n]  = C[n] + delta * C[n - 1];\\
S[n]  = (2 * n + 1 + delta * x * n) * C[n - 1] - (n + 1) * x * C[n - 2];\\
A = 2 * (n + 1) * (2 * n + 1);\\
\#\\
\#\\
\#\#\#\#\#\#\#\#\#\#\#\#\#\#\#\#\#\#\#\#\#\#\#\#\#\#\#\#\#\#\#\#\#\#\#\#\#\#\#\#\#\#\#\#\#\#\\
\# three families of C1 rules, see section 3. \#\\
\#\#\#\#\#\#\#\#\#\#\#\#\#\#\#\#\#\#\#\#\#\#\#\#\#\#\#\#\#\#\#\#\#\#\#\#\#\#\#\#\#\#\#\#\#\#\\
\#\\
\# get the orthogonal Gegenbauer polynomials for the weight (1 - x \^{} 2) \^{} 2 with:\\
\# C[n] := simplify (GegenbauerC (n, 5 / 2, x));\\
\#\\
\# C1, periodic with 1 interval, odd degree 2 * n + 1\\
\# see section 3.1.\\
\#\\
\# first quadrature rule, a point at the interval end, see section 3.1.1.\\
R[n - 1] = C[n - 1];\\
S[n - 2] = C[n - 2];\\
A = 2 * n * (n + 1) * (n + 2) / 9;\\
w[1] = 16 * (2 * n \^{} 2 + 6 * n + 1) / (3 * n * (n + 1) * (n + 2) * (n + 3));\\
\#\\
\# second quadrature rule, see section 3.1.2.\\
delta    = sqrt (3 * (n \^{} 2 + 3 * n - 1) / n / (n + 3));\\
R[n]     =   (n - 1) * (2 * n \^{} 2 + 2 * n - 3) * C[n]\\
           - (n + 3) * (2 * n \^{} 2 + 6 * n + 7 - 2 * (2 * n + 3) * delta) * C[n - 2];\\
S[n + 1] =   n * (6 * n \^{} 2 + 6 * n - 3 + 2 * (2 * n + 1) * delta) * (1 + x \^{} 2) * C[n - 1]\\
           + (2 * n \^{} 2 + 6 * n + 1) * (- 4 * (2 * n + 1) * x * C[n - 2]\\
                                        + (n + 2) * (1 + x \^{} 2) * C[n - 3]);\\
A = 2 * (n - 1) * (n + 1) * (n + 2) * (2 * n + 1) * (2 * n + 3)\\
      * (2 * n \^{} 2 + 2 * n - 3) * (2 * n \^{} 2 + 6 * n + 1) / 9;\\
\#\\
\# C1, periodic with 2 intervals, even degree 2 * n, "1/2 rules"\\
\# see section 3.2.\\
\#\\
delta = sqrt (3 * n * (n + 2) * (n \^{} 2 + 2 * n - 2));\\
\# 1. interval with n points\\
R[n - 1] =   (n - 1) * (2 * n \^{} 2 + 2 * n - 3) * C[n - 1]\\
           + (2 * delta + 3 - n - 6 * n \^{} 2 - 2 * n \^{} 3) * C[n - 2];\\
S[n - 2] =   (3 * (n + 2) * (2 * n \^{} 2 - 1) - 2 * delta) * C[n - 2]\\
           + (n + 2) * (2 * n \^{} 2 + 2 * n - 3) * C[n - 3];\\
A = 2 * (n - 1) * n * (n + 1) * (n + 2) * (2 * n + 1) * (2 * n \^{} 2 + 2 * n - 3) \^{} 2 / 9;\\
w[1] = 8 * (2 * n \^{} 2 + 4 * n - 3) * (2 * n \^{} 4 + 8 * n \^{} 3 + 4 * n \^{} 2 - 8 * n - 3 - delta)\\
         / (3 * (n - 1) * n * (n + 2) * (n + 3) * (2 * n + n \^{} 2 - 2) * (n + 1) \^{} 2);\\
\# 2. interval with n - 1 points\\
\# n - 1 points of the first interval reflected, without the first point at - 1\\

\normalsize \normalfont

\section{Visualisation of the quadrature rules}

\subsection{The rules with a period of one interval}

\begin{figure}[ht]
$$
\begin{array}{cc}
\includegraphics[width=0.4\textwidth]{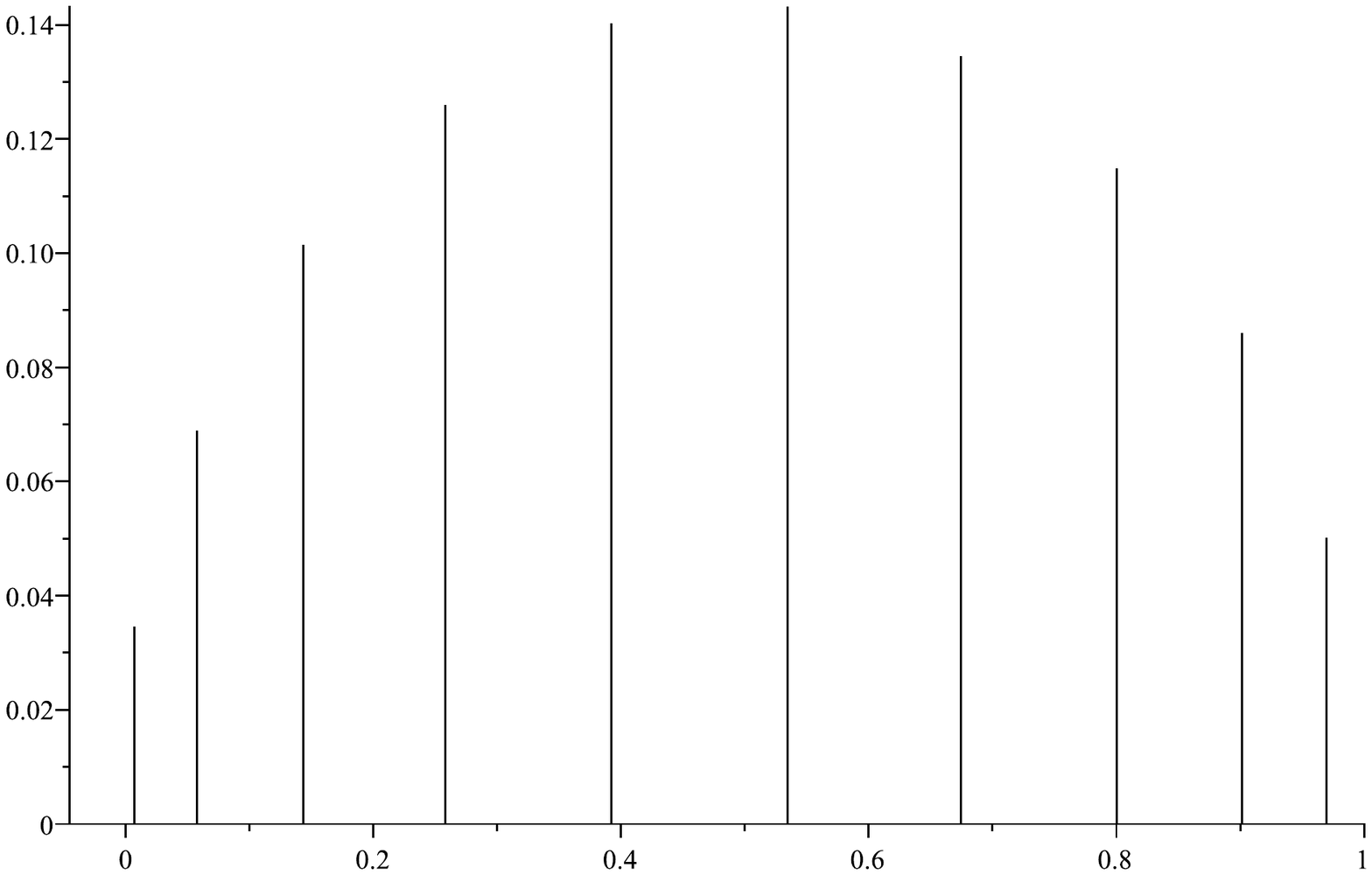} &      
\includegraphics[width=0.4\textwidth]{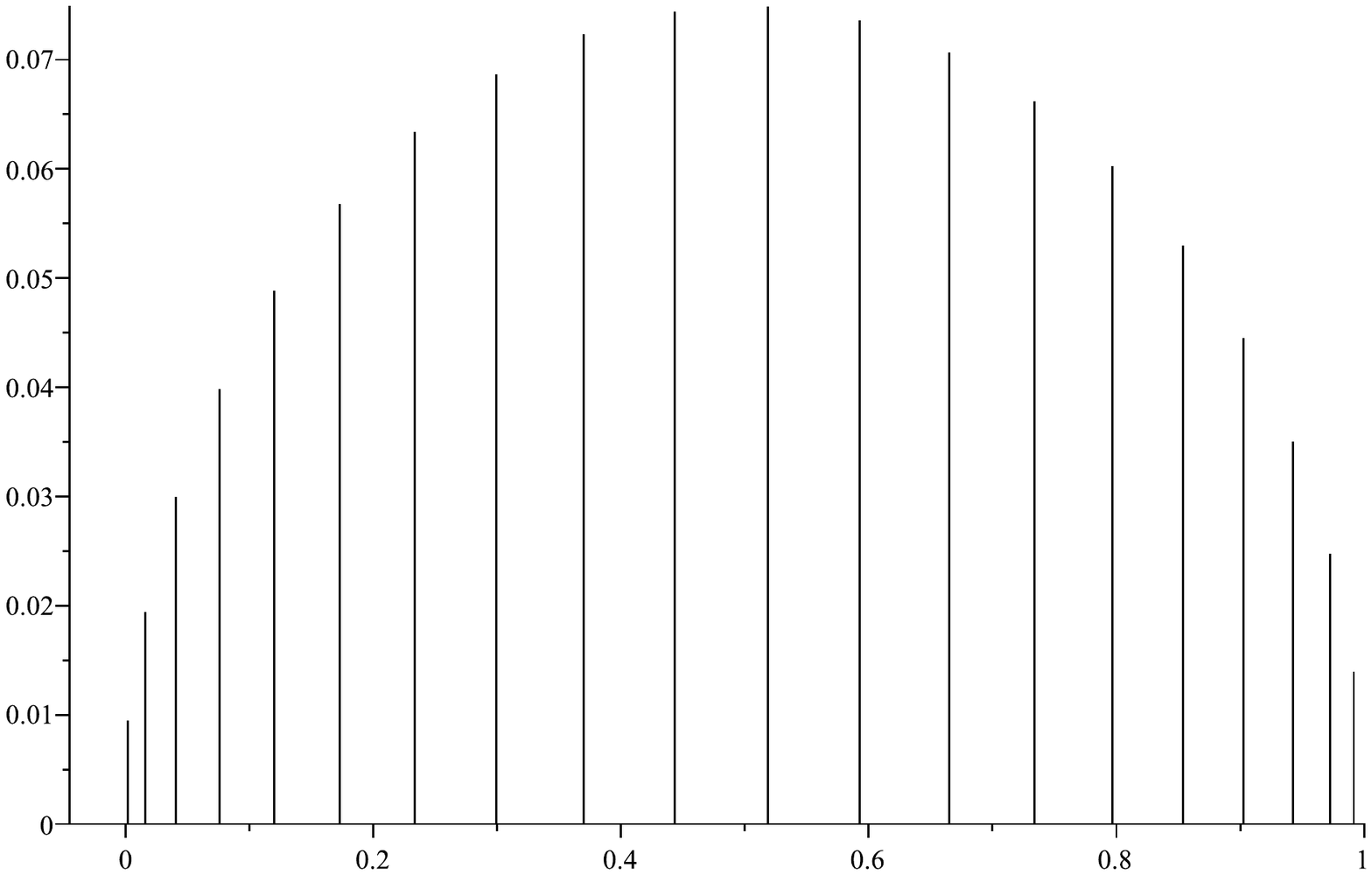} \\
(a)  &  (b) 
\end{array}
$$	
\caption{ c = 0, the weights of rule 2.2. 
$(a)$ rule with $n = 10$ points i.e. even degree $D = 2 n = 20$,
$(b)$ rule with $n = 20$ points i.e. even degree $D = 2 n = 40$, 
the rule has no reflection symmetry. }
\label{fig:C0xDEven1Int}
\end{figure}

\begin{figure}[ht]
$$
\begin{array}{cc}
\includegraphics[width=0.4\textwidth]{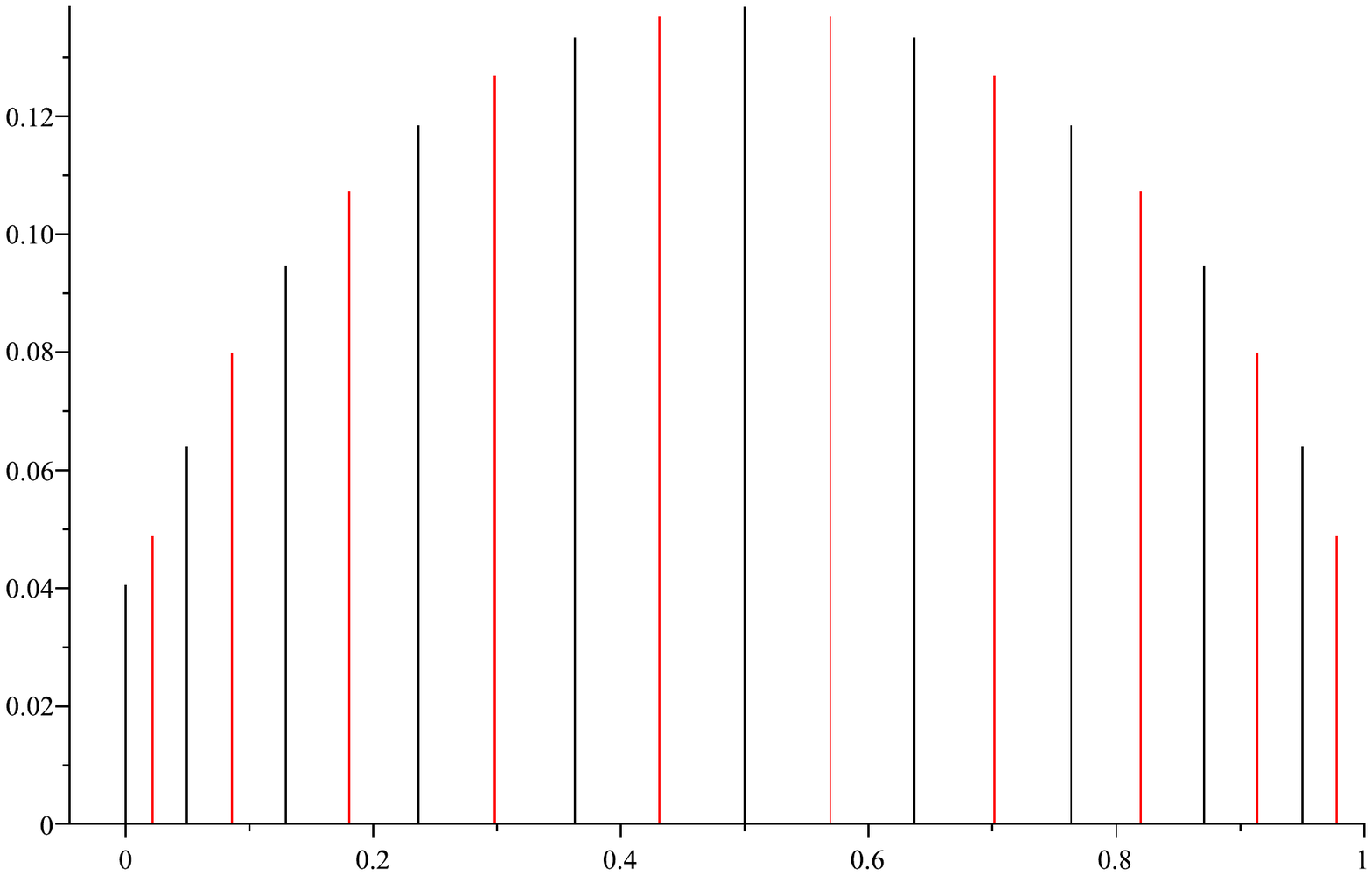} &      
\includegraphics[width=0.4\textwidth]{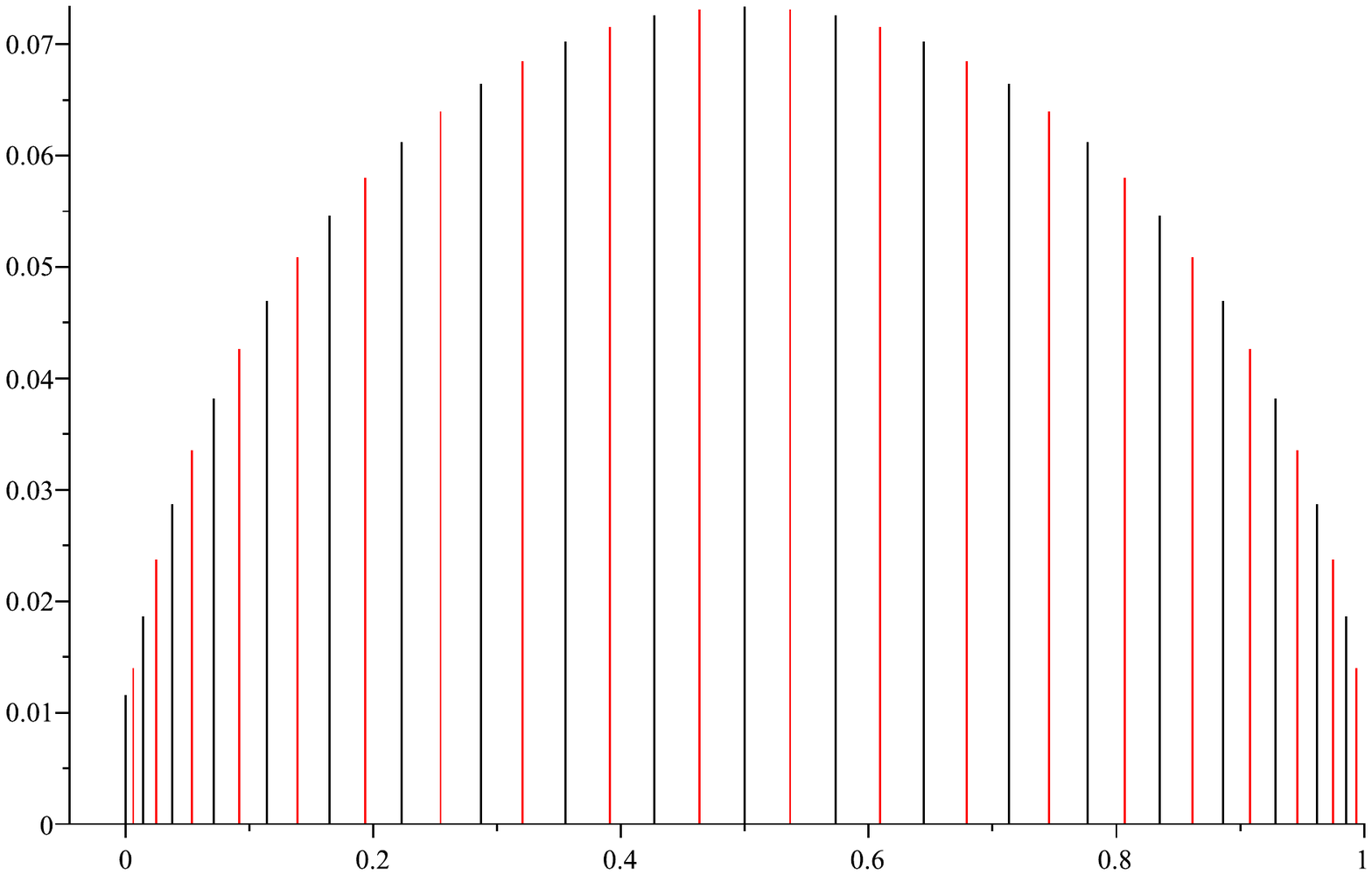} \\
(a)  &  (b) 
\end{array}
$$	
\caption{ c = 1, the weights of rule 3.1.1. in black, the weights of rule 3.1.2 in red
$(a)$ rule with $n = 10$ points i.e. odd degree $D = 2 n + 1 = 21$,
$(b)$ rule with $n = 20$ points i.e. odd degree $D = 41$, 
the rule has reflection symmetries at the interval boundaries and the midpoints. }
\label{fig:C1xDOddx1IntBoth}
\end{figure}

\newpage

\subsection{The rules with a period of two intervals}

\begin{figure}[ht]
$$
\begin{array}{c}
\includegraphics[width=0.8\textwidth]{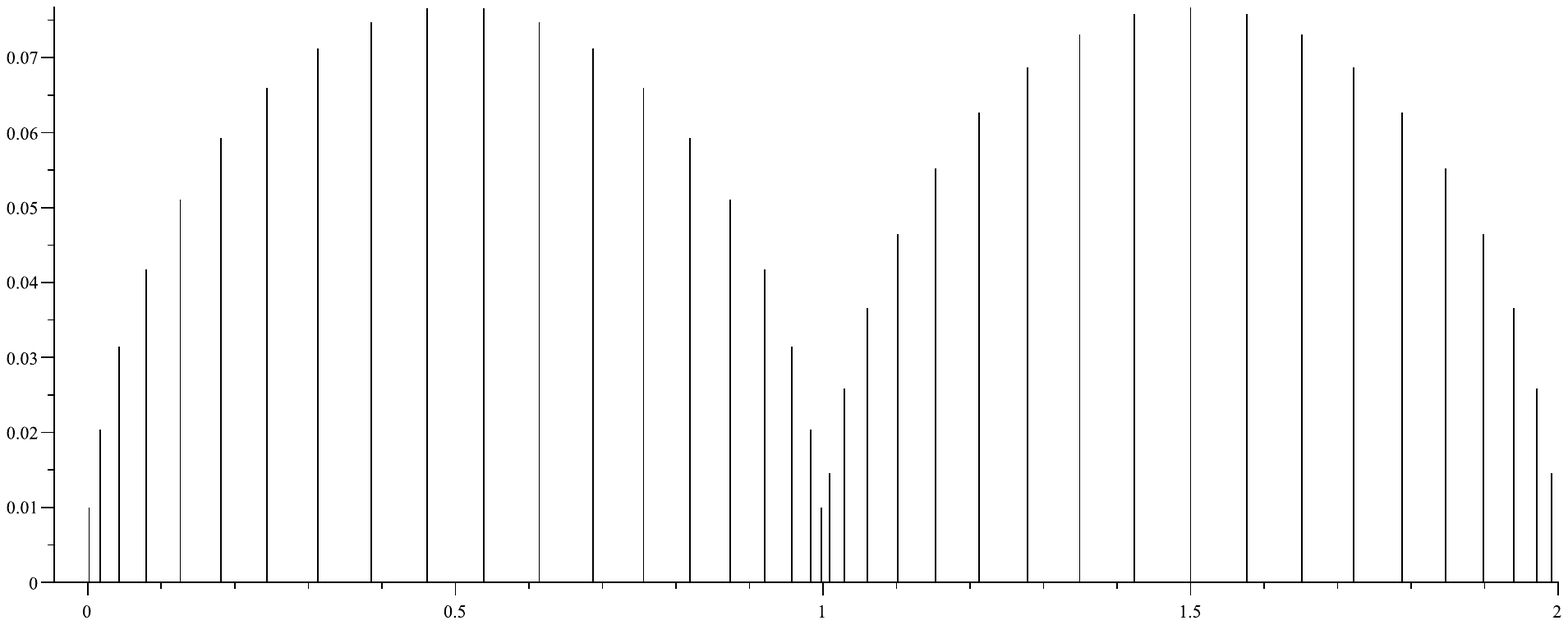}
\end{array}
$$	
\caption{ c = 0, the weights of rule 2.1. 
rule with $n = 20$ points in [0, 1] and 19 points in [1, 2] i.e. odd degree $D = 2 n - 1 = 39$, 
the rule has reflection symmetries at the interval midpoints. }
\label{fig:C0xDOdd2Int}
\end{figure}

\begin{figure}[ht]
$$
\begin{array}{c}
\includegraphics[width=0.8\textwidth]{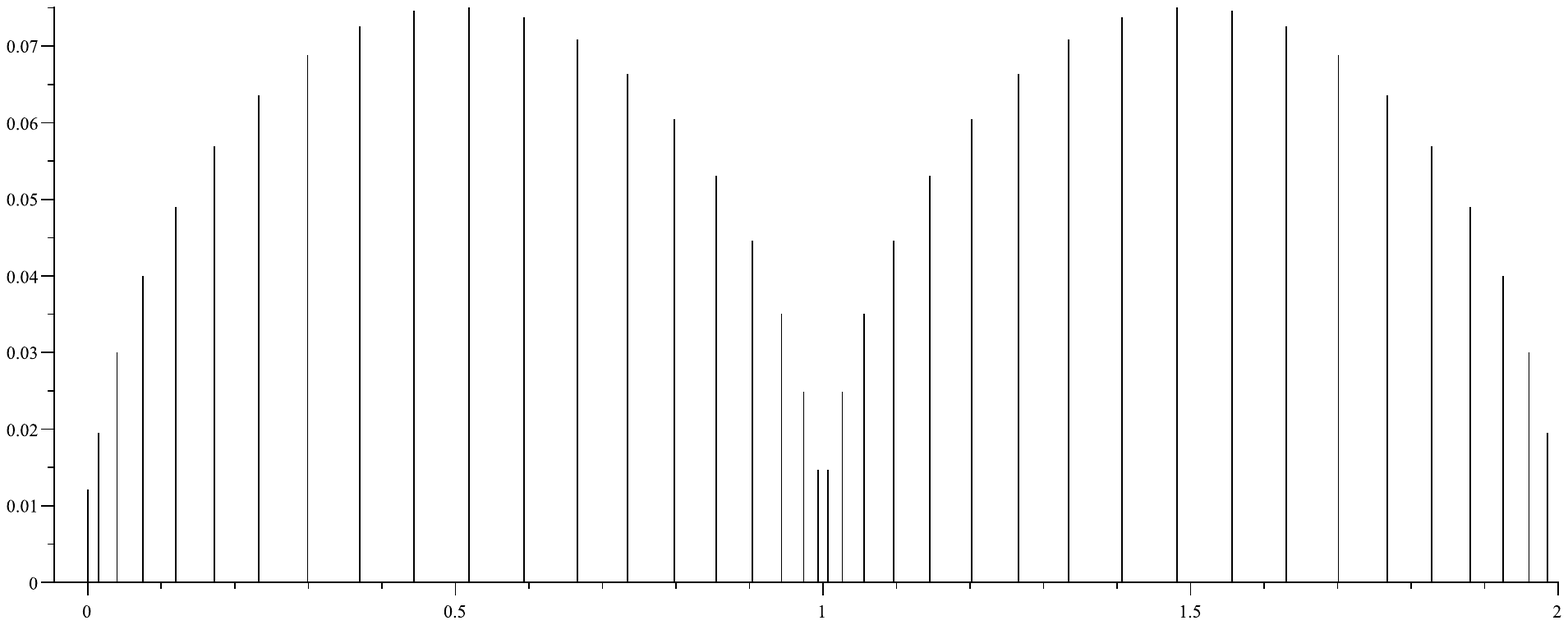}
\end{array}
$$	
\caption{ c = 1, the weights of rule 3.2.,
rule with $n = 20$ points in [0, 1] and 19 points in [1, 2] i.e. even degree $D = 2 n = 40$, 
the rule has reflection symmetries at the interval boundaries. }
\label{fig:C1xDEven2Int}
\end{figure}


\begin{thebibliography}{00}

\bibitem{BC2016}
Michael Barto\v{n}, Victor Manuel Calo, Gauss-Galerkin quadrature rules for quadratic and cubic spline spaces and their application to isogeometric analysis, Preprint, arXiv:1602.01200v1 [math.NA] 3 Feb 2016.
\bibitem{BC2015}
Michael Barto\v{n}, Victor Manuel Calo, Optimal rules for isogeometric analysis, Preprint, arXiv:1511.03882v1 [math.NA] 12 Nov 2015.
\bibitem{RU}
H. Ruhland, Quadrature rules for $C^{0}$ and $C^{1}$ splines, a recipe, unpublished manuscript June 2016.
\bibitem{HRS2010}
T.J.R. Hughes, A. Reali, and G. Sangalli, Effcient quadrature for NURBS-based isogeometric analysis, Computer Methods in Applied Mechanics and Engineering, 199 (58): 301-313, 2010.

\end{thebibliography}
\end{document}